\documentclass[11pt]{amsart}

\usepackage[T1]{fontenc}
\usepackage{ae,aecompl}

\usepackage{amssymb,amscd}
\usepackage{url}
\usepackage{amsthm,amsfonts,amsmath,amsxtra}
\usepackage{bbm}
\usepackage{graphicx}
\usepackage{epsfig,psfrag}

\numberwithin{equation}{section}
\numberwithin{table}{section}
\numberwithin{figure}{section}

\newtheoremstyle{bold}
{.5\baselineskip}{.5\baselineskip}{\itshape}{}{\bfseries}{.}{.5em}{}
\newtheoremstyle{shy}
{.5\baselineskip}{.5\baselineskip}{}{}{\bfseries}{.}{.5em}{}

\makeatletter
\def\@captionfont{\small}

\def\mychapter{%
  \if@openright\cleardoublepage\else\clearpage\fi
 \thispagestyle{empty}\global\@topnum\z@
  \@afterindenttrue \secdef\@mychapter\@schapter}

\def\@mychapter[#1]#2#3{\refstepcounter{chapter}%
  \ifnum\c@secnumdepth<\z@ \let\@secnumber\@empty
  \else \let\@secnumber\thechapter \fi
  \typeout{\chaptername\space\@secnumber}%
  \def\@toclevel{0}%
  \ifx\chaptername\appendixname \@tocwriteb\tocappendix{chapter}{#2\\ \scshape #3}%
  \else \@tocwriteb\tocchapter{chapter}{#2\\ \scshape #3}\fi
  \chaptermark{#1}%
  \addtocontents{lof}{\protect\addvspace{10\p@}}%
  \addtocontents{lot}{\protect\addvspace{10\p@}}%
  \@mymakechapterhead{#2}{#3}\@afterheading}

\def\@mymakechapterhead#1#2{\global\topskip 7.5pc\relax
  \begingroup
  \fontsize{\@xivpt}{18}\bfseries\centering
    \ifnum\c@secnumdepth>\m@ne
      \leavevmode \hskip-\leftskip
      \rlap{\vbox to\z@{\vss
          \centerline{\normalsize\mdseries
              \uppercase\@xp{\chaptername\ \thechapter}}
          \vskip 3pc}}\hskip\leftskip\fi
     #1\par \vskip 1pc
     \Large\mdseries\scshape\centering
     #2\par \endgroup
  \skip@34\p@ \advance\skip@-\normalbaselineskip
  \vskip\skip@ }

\def\section{\@startsection{section}{1}%
  \z@{.9\linespacing\@plus\linespacing}{.5\linespacing}%
  {\large\bfseries\boldmath\centering}}

\def\subsection{\@startsection{subsection}{2}%
  \z@{.7\linespacing\@plus\linespacing}{.5\linespacing}%
  {\normalfont\scshape\centering}}

\def\theindex{\@restonecoltrue\if@twocolumn\@restonecolfalse\fi
  \columnseprule\z@ \columnsep 35\p@
  \@indextitlestyle
  \thispagestyle{empty}%
  \let\item\@idxitem
  \parindent\z@  \parskip\z@\@plus.3\p@\relax
  \raggedright
  \hyphenpenalty\@M
  \footnotesize}

\renewcommand{\@bibtitlestyle}{%
  \@xp\section\@xp*\@xp{\bibname}%
}
\renewcommand{\tocchapter}[3]{%
  \indentlabel{\@ifnotempty{#2}{\ignorespaces#1 #2.\quad}}#3}

\renewcommand{\tocsection}[3]{%
  \indentlabel{\@ifnotempty{#2}{\makebox[3.2em][l]{\ignorespaces#1 #2.}}}#3}

\makeatother

\renewcommand{\tocappendix}[3]{%
  \indentlabel{#1.\quad}#3}
\renewcommand{\tocappendix}[3]{%
  \indentlabel{\makebox[5.7em][l]{\ignorespaces#1.}}#3}

\renewcommand{\bibname}{References}

\renewcommand{\geq}{\geqslant}

\renewcommand{\leq}{\leqslant}

\theoremstyle{bold}
\newtheorem{theorem}{Theorem}[section]
\newtheorem{proposition}[theorem]{Proposition}

\newtheorem{corollary}[theorem]{Corollary}

\theoremstyle{shy}

\newtheorem{remark}[theorem]{Remark}

\newcommand{\cA}{\mathcal{A}}
\newcommand{\cB}{\mathcal{B}}
\newcommand{\cC}{\mathcal{C}}

\newcommand{\cP}{\mathcal{P}}

\newcommand{\one}{\mathbbm{1}}

\newcommand{\EE}{\mathbb{E}}

\newcommand{\NN}{\mathbb{N}}

\newcommand{\PP}{\mathbb{P}\ts}

\newcommand{\dd}{\ts\mathrm{d}\ts}
\newcommand{\ee}{\ts\mathrm{e}\ts}

\newcommand{\ts}{\hspace{0.5pt}}
\newcommand{\nts}{\hspace{-0.5pt}}

\newcommand{\udo}[1]{\underaccent{$\text{.}$}{#1\ts}\nts}
\usepackage{accents}

\makeindex

\begin{document}

\title{Ancestral lines under recombination}

\author{Ellen Baake and Michael Baake}

\address{\{Faculty of Technology, Faculty of Mathematics\}, Bielefeld University, 
\hspace*{\parindent}Postbox 100131, 33501 Bielefeld, Germany}
\email{\{ebaake,mbaake\}@math.uni-bielefeld.de}

\maketitle

Solving the  recombination equation has been a long-standing challenge of \emph{deterministic} population genetics. We review recent progress obtained by introducing ancestral processes, as traditionally used in the context of \emph{stochastic} models of population genetics, into the deterministic setting. With the help of an ancestral partitioning process, which is obtained by letting population size tend to infinity (without rescaling parameters or time) in  an ancestral recombination graph, we obtain the solution to the recombination equation in a transparent form. 

\section[Introduction]{Introduction}
\label{sec:introduction}
\index{recombination} \emph{Recombination} is a genetic mechanism that `mixes' or `reshuffles' the genetic material of different individuals from generation to generation; it takes place in the course of  sexual reproduction. Models that describe the evolution of  populations
under recombination (in isolation or in combination with other processes) are among the major challenges in population genetics. Besides being of theoretical and mathematical interest, they play a major role in inference from population sequence data; compare the contribution of Dutheil \cite{EBMB-Du20} in this volume.

In line with the general situation in population genetics, models of recombination  come in two categories, deterministic and stochastic. In addition, there are versions in discrete and in continuous time, both of which will be considered below. In particular, our approach will result in a unified treatment of both.

Deterministic approaches assume that the population is so large that a law of large numbers applies and random fluctuations may be neglected. The resulting models are (systems of) \index{ordinary differential equation} ordinary   differential equations or (discrete-time) \index{dynamical system} dynamical systems, which describe the evolution of the genetic composition of a population under recombination in the usual forward direction of time; for review, see  \cite{EBMB-Ly92,EBMB-Ch99,EBMB-Bu00}. The genetic composition
is described via a probability distribution (or measure) on a space of
sequences of finite length.   The equations are nonlinear and notoriously difficult
to solve. Elucidating the underlying structure and finding solutions
was a challenge to theoretical population geneticists for nearly
a century. Indeed, the first studies go back to Jennings \cite{EBMB-Je1917} in 1917  and
Robbins \cite{EBMB-Ro1918} in 1918.  Geiringer \cite{EBMB-Ge44} in 1944  and
Bennett \cite{EBMB-Be54} in 1954  were the first to state the generic general
\emph{form} of the solution in terms of a convex combination of
certain basis functions, and evaluated  the corresponding coefficients recursively for sequences with a small number of sites.  The approach was later continued
within the systematic framework of genetic algebras; compare
\cite{EBMB-Ly92,EBMB-HaRi83}. The recursions for the coefficients were worked out in fairly general form by Dawson \cite{EBMB-Da02}. In any case, however, the work is technically cumbersome and yields limited insight into the underlying mathematical structure.

Stochastic approaches, on the other hand, take into account the fluctuations due to finite population size. The evolution of the composition of a population over time is described via a \index{Moran model!with recombination} Moran or a \index{Wright--Fisher model!with recombination} Wright--Fisher model with recombination. The first study goes back to Ohta and Kimura \cite{EBMB-OhKi69} in 1969.
Over the decades, two major lines of research have emerged. There has been continuous interest in  how the correlations between sites (known as  \index{linkage disequilibrium} linkage disequilibria) will develop; see \cite{EBMB-OhKi69,EBMB-SoSo07} and the overviews in
\cite[Ch.~5.4]{EBMB-HSW05}, \cite[Ch.~3.3 and 8.2]{EBMB-Du08} or \cite[Ch.~7.2.4]{EBMB-Wa09}.  The explicit time course of the genetic composition of the population is even more challenging, due to an intricate interplay of resampling
and recombination; compare~ \cite{EBMB-OhKi69,EBMB-SoSo07,EBMB-BaHe08,EBMB-BoKi10} as well as \cite[Ch.\ 8.2]{EBMB-Du08}.  These questions are usually approached forward in time.

The second line of research is concerned with genealogical aspects. Here, one starts with a sample taken from the present population and traces back the ancestry of the various sequence segments the individuals are composed of. The standard tool for this purpose  is the  \index{ancestral recombination graph} \emph{ancestral recombination graph} (ARG), first formulated by Hudson \cite{EBMB-Hu83} in 1983. 
His original
version was for two sites, but  generalisations to an arbitrary number of sites
\cite{EBMB-GrMa96,EBMB-BhSo12}  
and continuous versions 
\cite[Ch.~3.4]{EBMB-Du08} are immediate.

The stochastic models of recombination are related to their deterministic counterparts via a \index{law of large numbers!dynamical} dynamical law of large numbers as population size tends to infinity. Nevertheless, deterministic and stochastic approaches have largely led separate lives for  decades. It is the goal of this article to review recent progress  to build  bridges between them by introducing the genealogical picture into the deterministic equations. The corresponding ancestral processes remain random even in the deterministic limit, since they describe the history of single individuals (or a finite sample thereof). They lead to stochastic representations of the solutions of the deterministic equations and shed new light  both on their dynamics and their asymptotic behaviour.
A similar programme has been carried out for mutation-selection models; see \cite{EBMB-BCH18,EBMB-Co17} as well as the review \cite{EBMB-BaWa18}.

\section[Moran model with recombination]{Moran model with recombination}
\index{Moran model!with recombination}
Let us start from the \emph{Moran model with recombination} (in \emph{continuous time}),
which we recapitulate here from \cite{EBMB-BoKi10,EBMB-EPB16,EBMB-Es17}.
We consider a linear arrangement (or \index{sequence} \emph{sequence})
of $n$ discrete positions called \textit{sites}, which are collected in the set $S=\{1,\dotsc,n\}$. A site may be understood as a nucleotide site or a gene locus. We will throughout
consider sequences as (haploid) \emph{individuals}, that is, we think at
the level of gametes (rather than that of  diploid individuals that
carry two copies of the genetic information). 
Site $i$ is occupied by a letter $x_i \in X_i$, where  $X_i$ is a finite set, $1\leqslant i \leqslant n$. If sites are nucleotide sites, a natural
choice for each $X_i$ is the nucleotide alphabet $\{\rm{A,G,C,T}\}$; 
if sites are gene loci,
$X_i$ is the set of alleles that can occur at locus $i$.
The genetic type of each individual is thus described by the sequence $ x=(x_1,x_2,\dotsc,x_n) \in X_1 \times \dots \times X_n =\mathrel{\mathop:} X$, where
$X$ is the type space\footnote{We restrict ourselves to a finite type space here for ease of presentation; but the results generalise to more general type spaces where the $X_i$ may be locally compact \cite{EBMB-BaBa16}.}. 

\begin{figure}
\includegraphics[width=0.95\textwidth]{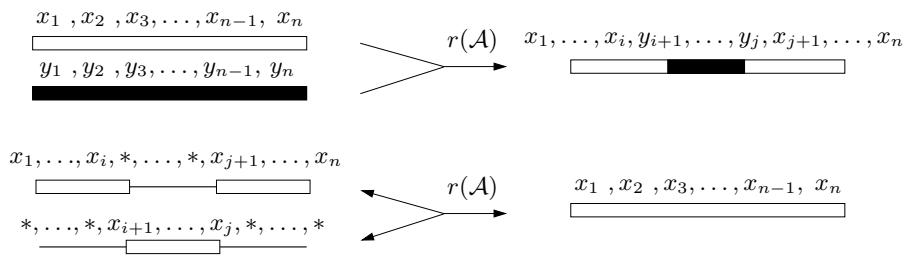}
\caption{\label{EBMB-recombining_sequences} Result of a double crossover
 between sites $i$ and $i+1$ and between $j$ and $j+1$ 
($1 \leqslant i<j < n$). Top: full details of parental sequences; bottom: a version that marginalises over the letters that do not end up in the offspring.
}
\end{figure}

In this setting, recombination means that a new individual
is formed as a `mixture' of an (ordered) pair\footnote{\label{EBMB-2parents}We formulate the model and the results   throughout for the (biologically realistic) case  of two parents here. Everything generalises easily to the situation with an arbitrary number of parents, which is mathematically interesting as well. Indeed, most of the results are available in the general setting in the original articles.} of parents, say $x$ and $y$.
We describe this mixture with the help of a \index{partition} \emph{partition} $\cA$ of $S$ into two parts. Namely, $\cA = \{A_1 , A_2\}$ means that the new individual copies the letters at all sites in $A_1$ from the first individual
and the letters at all sites in $A_2$ from the second; this happens via a number of \index{crossover} \emph{crossovers} between the sequences,  as illustrated in Figure~\ref{EBMB-recombining_sequences}. For reasons of symmetry, we need not keep track
of which part (or block) was `maternal' and which was `paternal'. 
Altogether, whenever an offspring is created, its sites are partitioned between parents according
to $\cA$ with probability $r(\cA)$, where $r(\cA) \geqslant 0$,
\mbox{$\sum_{\cA \in \cP_2(S)}r(\cA) \leqslant 1$}, and $\cP_2(S)$ is the set of all 
partitions of $S$ into two parts.
The sum $\sum_{\cA \in \cP_2(S)}r(\cA)$ is the probability that some recombination
event takes place during reproduction. With probability
$r(\boldsymbol{1})=1- \sum_{\! \cA \in \cP_2(S)}r(\cA)$, there is no recombination, 
in which case the offspring 
is the full copy of a single parent; here $\boldsymbol{1}\mathrel{\mathop:}= \{S\}$, the coarsest \index{partition} partition. 
We write $\cP_{\! \leqslant 2}(S)\mathrel{\mathop:}= \cP_2(S) \cup \{\boldsymbol{1}\}$ for the set of 
partitions into at most two parts, and $\cP(S)$ for the set of all partitions of $S$. The collection $\{r(\cA)\}_{\cA \in \cP_{\leqslant 2}(S)}$ is known as
the \index{recombination!distribution} \emph{recombination distribution}  \cite[p. 55]{EBMB-Bu00}.

\begin{figure}
\includegraphics[width=0.95\textwidth]{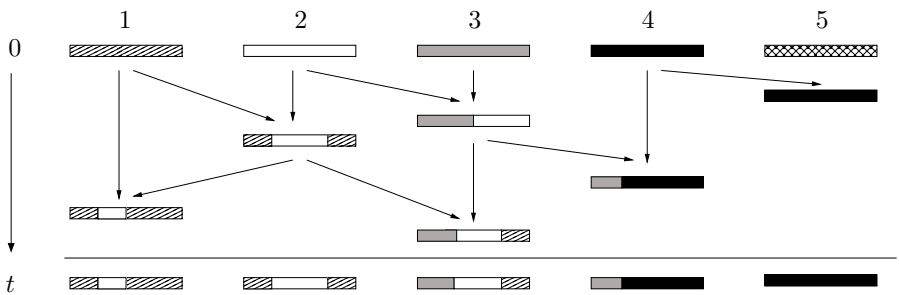}
\caption{\label{EBMB-moran} \index{Moran model!with recombination} A realisation of the Moran model with recombination forward in
time, with $N = 5$. For example, in the second event, individual 3  is replaced by a
recombined copy of individuals 2 and 3.}
\end{figure}

Consider now a \emph{population} of a constant number  $N$ of haploid individuals (that is, gametes) that evolves in continuous time as described next (see Figure~\ref{EBMB-moran}). Each individual dies at rate $\mu$, that is, it has an exponential lifespan with parameter $\mu$
(this parameter  simply sets the time scale).
When an individual dies, it is replaced by a new one as follows.
First, draw a \index{partition} partition $\cA$ according to the recombination distribution. 
Then, draw $|\cA|$ parents from the population (the parents may
include the  individual that 
is about to die), uniformly and with replacement, where $|\cA|$ is the number
of parts in $\cA$. If $|\cA|=2$, $\cA$ is of the form $\{A_1,A_2\}$, and the offspring inherits the sites
in $A_1$ from the first parent and the sites in $A_2$
from the second, as
described above.  If $|\cA|=1$ (and thus $\cA=\boldsymbol{1}$), 
the offspring is a full copy of a single parent 
(again chosen uniformly from among
all individuals); this
is called a (pure) \emph{resampling} event. All events
are independent of each other. Note that the model may equivalently be formulated in terms of 
reproducing rather than dying individuals, in the following way. Every individual reproduces at rate $\mu$, draws a partition $\cA$ according 
to the recombination distribution, and picks $|\cA|-1$ partners from the population; the offspring individual is pieced together
according to $\cA$ from the `active' individual and the partners, and replaces a uniformly chosen individual.

We identify the population at time $t$ with a  (random)  counting measure
$Z^{(N)}_t$ on $X$, where the upper index indicates the dependence on populaton size. Namely, $Z^{(N)}_t(\{x\})\geq 0$ denotes the number of individuals
 of type $x \in X$ at time $t$, and
$Z^{(N)}_t(\mathbb{A}) = \sum_{x \in \mathbb{A}} Z^{(N)}_t(\{x\})$ for $\mathbb{A} \subseteq X$. We can also write
\[
Z^{(N)}_t = \sum_{x \in X} Z^{(N)}_t(\{x\}) \, \delta_x
\]
 in terms of point measures on $x$. Since our Moran population has
constant size $N$, we have $\|Z^{(N)}_t\| = N$ for all times, where 
$\| Z^{(N)}_t \| \mathrel{\mathop:}= Z^{(N)}_t(X) = \sum_{x \in X} Z^{(N)}_t(\{x\})$
is the norm (or total variation) of $Z^{(N)}_t$.

This way, $(Z^{(N)}_t)^{}_{t \geqslant 0}$ is  a Markov chain in \emph{continuous time} 
with values in 
\begin{equation}\label{E}
E \mathrel{\mathop:}= \big\{z \in \{0, \dotsc, N\}^{\lvert X \rvert } : \| z \| = N \big\},
\end{equation}
where $\lvert X \rvert$ is the number of elements in $X$.
We will describe the action of recombination on (positive) measures with the help of so-called \index{recombinator} \emph{recombinators} as introduced in \cite{EBMB-BaBa03}. Let $\boldsymbol{M}_+(X)$ be the set of all positive, finite measures on $X$, where we understand $\boldsymbol M_+(X)$ to include the zero measure. Define the canonical \index{projection} projection  $\pi^{}_I \colon X \mapsto \mbox{\LARGE $\times$}_{i\in I} X_i =\mathrel{\mathop:} X^{}_I$, for $\varnothing \neq I\subseteq S = \{1, \ldots, n\}$, by $\pi^{}_I(x) = (x_i)_{i\in I} =\mathrel{\mathop:} x^{}_I$ as usual. 
For $\omega \in \boldsymbol M_+(X)$, the shorthand  $\omega^I \mathrel{\mathop:}= \pi_I^{}. \ts\omega=\omega \circ \pi_I^{-1}$ indicates the marginal measure with respect to the sites in $I$, where $\pi_I^{-1}$ is the preimage of $\pi_I^{}$, and the operation $.$ (where the dot is on the line and should not be
confused with a multiplication sign) denotes the  
\emph{pushforward}
of $\omega$  w.r.t.\ $\pi^{}_I$. In terms of coordinates, 
the definition may be spelled out as
\[ 
   \omega^I (x^{}_I ) = (\omega \circ \pi_I^{-1}) (x^{}_I ) 
  = \omega \big(\{x\in X :  \pi^{}_I(x) = x^{}_I \} \big), \quad 
  x^{}_I \in X^{}_I.
\]
Note that  
$ \omega^S = \omega $. 

Consider now $\cA = \{A_1, \ldots, A_m\}
\in \cP(S)$ and 
$\omega \in \boldsymbol M_+(X)$, and define
the \emph{recombinator} as 
\begin{equation}\label{EBMB-recombinator_2}
R_{\cA}(\omega) \mathrel{\mathop:}= \frac{1}{\|\omega \|^{m -1}}\,
\bigotimes_{A \in \cA} \omega^A,
\end{equation}
where $\otimes$ indicates the  product measure and the definition extends consistently to  $R_{\cA}(0) =0$. Note that $R_{\boldsymbol {1}}(\omega) = \omega$. Clearly, $\|R_{\cA}(\omega)\| = \|\omega \|$
 for all \mbox{$\omega \in \boldsymbol M_+(X)$}.
In particular, $R_{\cA}$ turns $\omega \neq 0$ into the (normalised) product measure
of its marginals with respect to the blocks in $\cA$. 
If $Z_t=z$ is the current population, then $\frac{1}{\| z \| }R_{\cA}(z)  = \frac{1}{N} R_{\cA}(z) $ is the
type distribution  that results when a hypothetical individual is created by drawing marginal types
(as encoded by $\cA \in \cP(S)$)
from the current population, uniformly and with replacement.

We now  use the recombinators to reformulate the Moran model with recombination in a compact way. Namely, since all individuals die at rate
$\mu$, the population loses type-$y$ individuals at rate $\mu Z^{(N)}_t(\{y\})$. Each loss is replaced by a new individual, which is sampled uniformly from $\frac{1}{N} R_{\cA}^{} (Z^{(N)}_t)$ with
probability $r(\cA)$ for  $\cA \in \cP^{}_{\! \leqslant 2}(S)$. 
Therefore, when $Z^{(N)}_t = z$, the transition to $z + \delta_x - \delta_y$
occurs at rate
\begin{equation}\label{lambda} 
\lambda^{(N)}(z;\, y,\, x) \,  = \! \sum_{\cA \in \cP^{}_{\! \leqslant 2}(S)} \frac{1}{N} \varrho(\cA) \big( R_{\cA}^{} (z) \big )(\{x\}) \,  z(\{y\}),
\end{equation}
where $\varrho(\cA) = \mu \, r(\cA)$ is a recombination \emph{rate} (in line with the continuous-time model)\footnote{Note that the meaning of $\varrho(\cA)$ as a  recombination rate is best understood by recalling the equivalent formulation of the model where every individual reproduces at \emph{rate} $\mu$ and then picks partition $\cA$ with \emph{probability} $r(\cA)$.}.
The summand for $\cA=\boldsymbol{1}$ corresponds to pure resampling, whereas all
other summands include recombination.
Note that $\lambda^{(N)}$   includes `silent transitions' ($x=y$).

\textbf{Law of large numbers.}
Consider now the family of processes $(Z_t^{(N)})^{}_{t\geqslant 0}$ with
$N \in \NN$. Also, consider the normalised version 
$(\frac{1}{N} Z_t^{(N)})^{}_{t\geqslant 0}$; note that $\frac{1}{N} Z_t^{(N)}$ is a random probability
measure on $X$. For $N \to \infty$ and without any rescaling of the
$\varrho(\cA)$ or of time, the sequence $(\frac{1}{N} Z_t^{(N)})^{}_{t\geqslant 0}$ converges
to the solution of the \index{recombination!equation} \index{ordinary differential equation} \emph{deterministic recombination equation\footnote{\label{EBMB-generalpart}The generalisation to an arbitrary number of parents, that is $\cA \in \cP(S)$, is treated in \cite{EBMB-BaBa16}. The special case  $\cA \in \cP^{}_{\! \leq 2}(S)$ is then obtained by setting $\varrho(\cA)=0$ for all $\cA \notin \cP^{}_{\! \leq 2}(S)$. In any case, note that the summand for $\cA=\boldsymbol{1}$ may or may not be included in the right-hand side of the equation, since it contributes nothing due to $R_{\boldsymbol{1}}(\omega)=\omega$.}}
\begin{equation}\label{EBMB-detreco}
\dot\omega_t \, = \! \! \! \!  \sum_{\cA \in \cP^{}_{\! 2} (S)} \varrho(\cA) \big (R_{\cA}^{} (\omega_t) - \omega_t \big)
\end{equation}
with initial value $\omega^{}_0 \in \boldsymbol{P}(X)$ (the set of probability measures on $X$), where  we assume that
\[
  \lim_{N \to\infty}  \frac{Z_0^{(N)}}{N} \, = \, \omega^{}_0.
\]
The convergence to the differential equation~\eqref{EBMB-detreco}   is a 
\index{law of large numbers!dynamical} \emph{dynamical law of large numbers} and due to
\cite[Thm.~11.2.1]{EBMB-EtKu86}. 
The precise statement  as well as the proof are perfectly analogous to  \cite[Prop.~6]{EBMB-BaHe08}, 
which assumes a slightly different recombination and sampling scheme, without consequence for the convergence claim.

\section{Ancestral recombination graph and  deterministic limit}
Let us  return to the finite-$N$ model and construct the type of an individual sampled randomly   from the  population at time $t$ (the `present') by genealogical means. We do so by adapting the \index{ancestral recombination graph} ARG (see \cite{EBMB-BhSo12} and, for overviews,  \cite[Ch.~5.4]{EBMB-HSW05}, \cite[Ch.~3.3, 8.4]{EBMB-Du08} or \cite[Ch.~7.2.4]{EBMB-Wa09}) to our model and a sample of size 1.

\begin{figure}
\includegraphics[width=0.95 \textwidth]{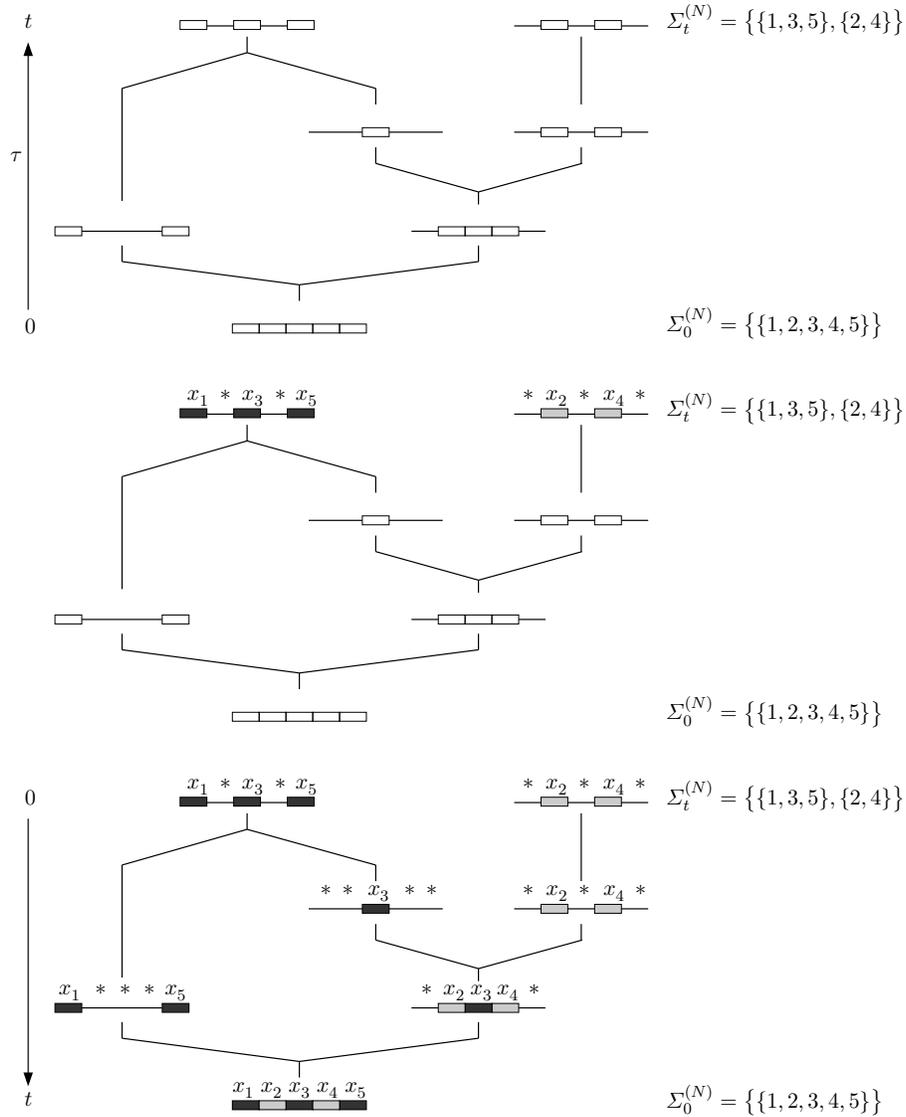} 
\caption{\label{EBMB-arg} Example realisation of the \index{partitioning process} partitioning process (top), assigning letters to the parts (middle), and propagating them downwards (bottom).}
\end{figure}

The type of an individual at present, together with its ancestry, can thus be constructed by a three-step procedure as illustrated in Figure~\ref{EBMB-arg}. 
First, we run a \index{partitioning process} partitioning process $(\varSigma^{(N)}_\tau)_{0 \leq \tau \leq t}$. Here,  $(\varSigma^{(N)}_\tau)_{\tau \geq 0}$ is a Markov chain in continuous time on $\cP(S)$, whose time axis is directed  into the past; we use the variables $t$ and $\tau$ throughout for forward and backward time, respectively, so $\tau = t$ in backward time corresponds to $t = 0$ in forward time. 
\index{partitioning process} The  process starts  with the coarsest partition $\varSigma^{(N)}_0=\boldsymbol{1}$, that is, we consider the (intact) sequence of one individual at time $t$. Then, $(\varSigma^{(N)}_\tau)_{\tau \geq 0}$ describes the partitioning of sites into parental individuals at time $\tau$ before the present; sites in the same block (in different blocks) belong to the same (to different)
individuals. Clearly, $\lvert \varSigma^{(N)}_\tau \rvert$ is
the number of ancestral individuals at time $\tau$. The process $(\varSigma^{(N)}_\tau)_{\tau \geqslant 0}$ consists of   \emph{splitting}  
and \index{coalescence} \emph{coalescence}  events (and combinations thereof), is independent of the types, and will be described in detail below.

In the second step, a letter is assigned to each site of $S$ at  $\tau=t$  (that is, at forward time 0) as follows. For every part of
$\varSigma^{(N)}_t$,  pick an individual from the initial population $Z^{(N)}_0$ (without replacement) and copy its letters to the sites in the block considered. For illustration, also assign  a colour to each block, thus indicating different parental individuals. In the last step, the letters and colours are propagated downwards (that is, forward in time) according to the realisation of  $(\varSigma^{(N)}_\tau)_{0\leq \tau\leq t}$ laid down in the first step.

\begin{figure}
\includegraphics[width=.6\textwidth]{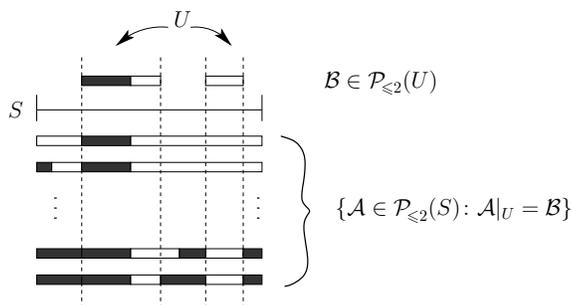}
\caption{\label{EBMB-fig:marginal_recorate} The marginal recombination rate
for a partition $\cB$ of a subset $U$ is the sum of all recombination rates
for partitions $\cA$ of the original set $S$ that lead to $\cB$ under
\index{projection} projection to $U$.}
\end{figure}

Let us now describe the \index{partitioning process} partitioning process more precisely, following \cite{EBMB-EPB16,EBMB-Es17} but specialising to $\varSigma_0^{(N)}=\boldsymbol{1}$.
Since we also trace back sites in subsets 
$U \subseteq S$ (rather than complete sequences), we need the corresponding \index{recombination!rates, marginal} \emph{marginal recombination rates}
\begin{equation}\label{EBMB-marg-rates}
  \varrho^U\! (\cB) \, \mathrel{\mathop:}= \!
  \sum_{\substack{\cA \in \cP^{}_{\! \leqslant 2}(S) \\ \cA |^{}_{U} = \cB}} 
  \! \varrho(\cA)
\end{equation}
for any $\cB\in\cP^{}_{\! \leqslant 2}(U)$,  where $\cA |^{}_{U}$ is the partition of $U$ induced by $\cA$; namely,   $\cA |^{}_{U} = \{ A \cap U : A \in \cA, A \cap U \neq \varnothing \}$.
Clearly, $\varrho^S(\cB)= \varrho(\cB)$ and $\varrho^U\!(\cB)$ is the sum of the rates of all recombination events that 
lead to partition $\cB$ under the \index{projection} projection to $U$, as illustrated in 
Figure~\ref{EBMB-fig:marginal_recorate}. Note that, for $|U|=1$, the only recombination parameter is $\varrho^U\!(\boldsymbol{1})=1$ (note that  we use  $\boldsymbol{1}$ to indicate the coarsest partition throughout, where the meaning is always clear from the upper index, so here $\boldsymbol{1}= \{U\}$).

\index{partitioning process} Suppose now that the current state is $\varSigma^{(N)}_\tau = \cA=\{A_1,\dotsc,A_m\}$ and denote by $\Delta$ the waiting time to the next event. The random variable $\Delta$ is exponentially distributed with parameter $m \mu$, since each block corresponds to an individual, and each individual is independently affected at rate $\mu$. When the event happens, choose one of the $m$ blocks, each with probability $\frac{1}{m}$.
If $A_j$ is picked,  $\varSigma^{(N)}_{\tau+\Delta}$ is obtained via a two-step procedure, namely a splitting step followed by a sampling step, namely: 
\begin{enumerate}
  \item In the \emph{splitting step}, block $A_j$ turns into an intermediate state $\mathfrak{a}$
with probability $r^{A_j}(\mathfrak{a})$, $\mathfrak{a} \in \cP^{}_{\! \leqslant 2}(A_j) $, where the marginal probabilities $r^U(\cB)$  are defined as the marginal recombination rates in Eq.~\eqref{EBMB-marg-rates} with $\varrho$ replaced by $r$. In detail: 
\begin{itemize}
\item \label{split_into_one} With probability $r^{A_j}(\boldsymbol{1})$, the block
$A_j$ remains unchanged. The resulting intermediate state (of this block) is $\mathfrak{a}=\boldsymbol{1}|_{A_j} = \{A_j\}$. 
\item \label{split_into_two} With probability $r^{A_j}(\mathfrak{a})$, $\mathfrak{a}\in \cP^{}_{2}(A_j)$,
 block $A_j$ splits into two parts, $\mathfrak{a}=\{A_{j_1},A_{j_2}\}$. 
\end{itemize}
\item In the following \emph{sampling step}, each block of $\mathfrak{a}$ chooses out of $N$ parents, uniformly and with replacement. Among these, there are $m-1$ parents that carry one block of $\cA \setminus \{A_j\}$ each; the remaining $N-(m-1)$ parents are \emph{empty}, that is, they do not carry  ancestral material available for coalescence. Coalescence happens if the choosing block  picks a parent that carries ancestral material; otherwise, the choosing block  becomes an ancestral block of its own, which is available for coalescence from then onwards. The possible outcomes are certain coarsenings of $(\cA \setminus \{A_j\}) \cup \mathfrak{a}$.
\end{enumerate}

The long list of outcomes is provided explicitly in \cite{EBMB-EPB16} for the special case of single crossovers and in \cite{EBMB-Es17} for general partitions into two parts, and the formal \index{duality} duality between the Moran model and the partitioning process is established. Here,   we only aim at the law of large numbers, which is again obtained as  $N \to \infty$ without rescaling of parameters or time. In this limit, each of the blocks of the intermediate state $\mathfrak{a}$ ends up in a different individual, so there are no coalescence events and $\mathfrak{a}$ is the final state. As a consequence, the blocks of the partition are conditionally independent.
This leads to the following result.

\begin{proposition}[Law of large numbers for the ARG \cite{EBMB-EPB16,EBMB-Es17}]\label{EBMB-detlimit}
\index{law of large numbers} \index{ancestral recombination graph}
The sequence of \index{partitioning process} partitioning processes 
$(\varSigma^{(N)}_\tau)_{\tau\geqslant 0}$, with $N \in \NN$ and initial state $\varSigma^{(N)}_0 \equiv \boldsymbol{1}$, converges in distribution, as $N \to \infty$,  to the
process $(\varSigma_\tau)_{\tau\geqslant 0}$ with initial state 
$\varSigma_0 = \boldsymbol{1}$
and generator matrix \/ $Q$ defined
by the nondiagonal elements
\[
 Q^{}_{\cA\cB} = \begin{cases} \varrho^{A}(\mathfrak{a}), & \text{if }  
                   \cB = (\cA \setminus \{A\})\cup \mathfrak{a} \text{ for some }
                   A \in \cA \text{ and } \mathfrak{a} \in \cP^{}_{\! 2} (A),
                   \\[2mm]
      0, &  \text{for all other } \cB \neq \cA.
\end{cases} 
\]
The limiting process may therefore be described as follows. If the current state is
$\varSigma_\tau=\cA$,  each part $A$ of $\cA$ is replaced by $ \mathfrak{a} \in \cP(A) \setminus \{ A \}$
at rate $\varrho^{A} (\mathfrak{a})$, independently of all other
parts.
Hence, $(\varSigma_\tau)_{\tau\geqslant 0}$ is a process
of  progressive \index{refinement} refinements, that is,
$\varSigma_T\preccurlyeq \varSigma_\tau$ for all $T \geq \tau$.    
\end{proposition}

\begin{figure}[t]
\index{partitioning process}
\includegraphics[width=0.95\textwidth]{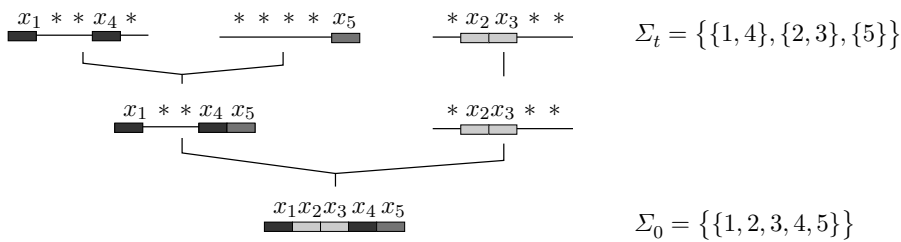}
\caption{\label{EBMB-partitioning_process} Determining the type of an individual at time $t$ via the \index{partitioning process} partitioning process $(\varSigma_\tau)_{0 \leq \tau \leq t}$.}
\end{figure}

\index{partitioning process} The process $(\varSigma_\tau)_{\tau\geqslant 0}$, which is  illustrated in Figure~\ref{EBMB-partitioning_process}, may be understood 
as the $N \to \infty$ limit of  the ARG started with a single individual. Note that, due to the continuous-time setting, at most one block may be refined at any given time (with probability one), but it may be any of the blocks.

Since $Q$ is the Markov generator of
$(\varSigma_\tau)_{\tau\geqslant 0}$, we can further conclude that
\[
    (\ee^{\tau \ts Q})^{}_{\cB \cC} \, = 
    \, \PP \bigl( \varSigma_{\tau} = 
    \cC \mid \varSigma_0 = \cB \ts \bigr )
\]    
(where $\PP$ denotes probability), that is, the transition
probability from  $\cB$ to  $\cC$ during a time interval
of length $\tau$. This leads us to the solution of the deterministic recombination equation.

\begin{theorem}[Solution of the recombination equation \cite{EBMB-BaBa16}\footnote{In fact, \cite{EBMB-BaBa16} treats the general case of an arbitrary number of parents, which corresponds to allowing for multiple (rather than binary) splits in the partitioning process; compare Footnotes~\ref{EBMB-2parents} and~\ref{EBMB-generalpart}.}]
  \label{EBMB-thm-recosol}
  The solution of the recombination equation \eqref{EBMB-detreco} reads 
\[
     \omega^{}_{t} \, =  
     \! \sum_{\cA\in\ts\cP(S)}\! \!  a^{}_t(\cA) 
    \, R^{}_{\cA} (\omega^{}_{0}) \, = \, \EE \big ( R_{\varSigma_t}(\omega_0)   \mid \varSigma_0= \boldsymbol{1} \big ),
\]
  where
  \[ a_t(\cA) \, = \, \PP \bigl( \varSigma_{t} = 
    \cA \mid \varSigma_0 = \boldsymbol{1} \bigr) = (\ee^{t \ts Q})^{}_{\boldsymbol{1} \cA}
  \]
  and \/ $\EE$ denotes expectation. 
\end{theorem}
\begin{remark}
With Theorem~\ref{EBMB-thm-recosol}, we have found a 
\emph{stochastic representation} of the  solution of the (deterministic) differential equation
\eqref{EBMB-detreco}. This reflects the fact that, while the time evolution
of the composition of the infinite population follows a (dynamical)
law of large numbers and is hence deterministic, the fate and ancestry
of a single individual retains some stochasticity. While
ancestral processes are common tools when working with the stochastic
processes that describe finite populations, they are not within the usual
scope of deterministic population genetics. 
\end{remark}

\begin{remark}
In \cite{EBMB-BaBa16}, the route of thought was, in fact, different from the one presented here. While we start from the ancestral process in this review, \cite{EBMB-BaBa16} works forward in time by means of classical methods from the theory of differential equations. The key was to establish the system of differential equations for the quantities $R_\cA (\omega_t)$, $\cA \in \cP(S)$, by exploiting the properties of the recombinators. This procedure mimics the algebraic technique of \index{Haldane linearisation} \emph{Haldane linearisation} and leads to the generator $Q$ in a purely analytic way. The partitioning process then emerged as an \emph{interpretation} of the result. 
\end{remark}

\begin{remark}
  It is easy to see that Theorem~\ref{EBMB-thm-recosol} extends to the 
  duality relation
  \[
    \EE \big ( R_\cB(\omega_t)   \mid \omega_0=\nu \big ) \, = \, 
    \EE \big ( R_{\varSigma_t}(\nu)   \mid \varSigma_0= \cB \ts \big )
  \]
  for any $\nu \in \boldsymbol{P} (X)$ and $\cB \in \cP(S)$. Hence, since the
  left-hand side is deterministic,
  \[
    R_\cB(\omega^{}_{t}) \, =  \,
    \EE \big ( R_{\varSigma_t}(\omega_0)   \mid \varSigma_0= \cB \ts \big ) 
  \]
  for any initial condition $\omega_0 \in \boldsymbol{P} (X)$.
 \end{remark}

Let us now turn  to the  evaluation of the $a_t$ of  
Theorem~\ref{EBMB-thm-recosol}.
It has been shown\footnote{The result in \cite{EBMB-BBS16} is again more general since it is not restricted to binary splitting.} in \cite{EBMB-BBS16} that,
in the generic case that the $\psi^{\ts U} \! (\cA)$ explained below
are all distinct, it can be given in the form
\begin{equation}\label{EBMB-gen-ansatz}
    a_{t}  (\cA) \, = \, \sum_{\udo{\cB}\succcurlyeq \cA} 
    \theta^{\ts S}\! (\cA,\cB\ts ) \, \ee^{- \psi^{\ts S}\! (\cB\ts )\ts  t}.
\end{equation}
Here, the underdot denotes the summation variable, $     \psi^{\ts U} \! (\boldsymbol{1}) \, \mathrel{\mathop:}= \, 
     \sum_{\cA\ne \boldsymbol{1}} \varrho^{U} \! (\cA) $
for all $\varnothing \neq U \subseteq S$, 
and the values for all other $\cA \in \cP(U)$ are  defined recursively by
$
   \psi^{\ts U} \! (\cA) \, \mathrel{\mathop:}= \sum_{i=1}^{|\cA|} \psi^{A_i} (\boldsymbol{1}) \ts .
$ 
In the context of the \index{partitioning process} partitioning process, $\psi^{A_i} (\boldsymbol{1})$ is 
the total rate of
any further partitioning of part $A_i$, and so, due to the
independence of the parts,  $\psi^{\ts U} \! (\cA)$
is the total rate of transitions out of state $\cA$.
The coefficients $\theta^{\ts U}\! (\cA,\cB\ts )$
follow the  recursion
\begin{equation}\label{EBMB-theta-rec}
   \theta^{\ts U} \! (\cA,\cB\ts ) \, = \! 
   \sum_{\cB\preccurlyeq\udo{\cC}\prec\boldsymbol{1}}
   \frac{\varrho^{U} \! (\cC)}{\psi^{\ts U}\! 
    (\boldsymbol{1}) - \psi^{\ts U} \! (\cB\ts )}
   \,\prod_{i=1}^{|\cC|} \theta^{C_i} \big ( \cA|^{}_{C_i} , \cB|^{}_{C_i} \big )
\end{equation}
for all $\cA\preccurlyeq \cB \prec \boldsymbol{1}$, where  the initial conditions are given by
$\theta^{\ts U} \! (\boldsymbol{1},\boldsymbol{1}) = 1$  
together with
$
    \theta^{\ts U} \! (\cA,\boldsymbol{1}) \, = \, - \! \!
    \sum_{\cA \preccurlyeq \udo{\cC}\prec \boldsymbol{1}}
    \! \theta^{\ts U} \! (\cA,\cC) \ts 
$
for $\cA \prec \boldsymbol{1}$ and all $U \subseteq S$.
Note that  everything is uniquely determined by the initial conditions for the singleton sets $U = \{i\}$ with $i \in S$.

The recursion exploits the lower-triangular form of $Q$. 
This  type of solution was motivated by earlier work of 
Geiringer \cite{EBMB-Ge44}, Bennett \cite{EBMB-Be54},
Lyubich \cite{EBMB-Ly92} and Dawson \cite{EBMB-Da02}, who worked on the analogous  system in
\emph{discrete} time (see below). We have made progress here by treating
the problem within a systematic
lattice-theoretic setting, which is the key for the
transparent construction of the solution. Furthermore, the measure-theoretic framework allows to also work with more general type spaces, where the $X_i$ may be locally compact \cite{EBMB-BaBa16}.

Let us note that the recursion \eqref{EBMB-theta-rec}
is of a fairly simple structure and computationally 
convenient.   In the next
section, we shall present an explicit solution for the special case of single-crossover recombination.

\begin{remark}
  Given the generator $Q$ from Proposition~\ref{EBMB-detlimit}, the matrix function of transition probabilities, $M(t) \mathrel{\mathop:}= \ee^{tQ}$, solves the Cauchy problem $\dot M = MQ$ with initial condition $M(0)=\one$ (where $\one$ is the identity matrix) and constitutes a Markov semigroup, so $M(t+s) = M(t) M(s)$ for $t,s \geq 0$. More generally, it is also of interest to consider the inhomogeneous counterpart, where $Q=Q(t)$ is time dependent; see \cite[Addendum]{EBMB-Ba05} for an example in the case of single crossovers. Let $M(t)$ again denote the solution of the Cauchy problem, which is unique under mild assumptions on $Q$ by general principles \cite{EBMB-EnNa00}. Clearly, $M(t)$ is still the matrix of transition probabilities until time $t$ and the underlying process satisfies the Markov property, while the semigroup property is lost.

  There are now two scenarios to be distinguished as follows. When the generator family $(Q(t))_{t \geq 0}$ is commuting, so \index{commutator} $Q(t)Q(s)=Q(s)Q(t)$ for all $t,s\geq 0$, one gets
\begin{equation}\label{EBMB-Mt}
  M(t) \, = \, \exp \int_0^t Q(\zeta) \dd \zeta
\end{equation}
or, more generally, $M(t,s) = \exp \int_t^s Q(\tau) \dd \tau$, with $M(t,s) M(s,r) = M(t,r)$ for $r \geq s \geq t \geq 0$, also known as the \index{flow} flow property. In general, however, the generators $Q(t)$ need not commute, and Eq.~\eqref{EBMB-Mt} has to be replaced by the more general \index{Peano-Baker} \emph{Peano--Baker formula}; see 
\cite{EBMB-BaSc11} for details. It can still be evaluated in some simple cases, and the flow property remains valid.
\end{remark}

Let us finally turn to the asymptotic behaviour of the solution of the recombination equation. It can, of course, be read off Eq.~\eqref{EBMB-gen-ansatz}, but it is more instructive to argue directly on the grounds of  $(\varSigma_t)_{t \geq 0}$. The following consequence of Proposition~\ref{EBMB-detlimit} and Theorem~\ref{EBMB-thm-recosol} is then immediate.

\begin{corollary}[Asymptotic behaviour of recombination equation]\label{EBMB-asymptotics}
Assume without loss of generality that $\varrho^{\{i,i+1\}}_{\{\{i\},\{i+1\}\}}>0$ for all $i \in S \setminus \{n\}$ $($if this is not the case, glue sites $i$ and $i+1$ together so that they form a single site$\ts)$. The \index{partitioning process} partitioning process is then absorbing, with
\[
 \lim_{t \to \infty}  \varSigma_t \, = \, \big \{ \{1\},\{2\}, \ldots, \{n\} \big \}
\]
almost surely and independently of $\varSigma_0$, and \index{equilibrium}
\[
  \lim_{t \to \infty} \omega_t^{}  \, = \, \bigotimes_{i=1}^{n} (\pi^{}_i . \omega^{}_0).
\]
\index{convergence!rate of} The convergence to the limit is exponentially fast.
\end{corollary}

\section{An explicit solution for single-crossover recombination}
\label{sec:single-cross}
There is an important special case that 
allows for
a \emph{closed} solution of the Markov semigroup, beyond the somewhat deceptive notation $\ee^{tQ}$ for the Markov semigroup generated by $Q$. This is the case of \emph{single crossovers}
of two parental gametes, which
is also highly relevant biologically: Since crossovers are rare,
it is  unlikely that two or more of them happen
in a given reproduction event, in any sequence of moderate length.

We speak of \index{recombination!single-crossover} \index{crossover!single} \emph{single-crossover recombination} if
$\varrho(\cA) > 0$ implies $\cA \! \in \!
\mathcal{I}_{\leqslant 2}(S)$. Here, $\mathcal{I}(S)$ is the set of
\emph{interval partitions} of $S$, $\mathcal{I}_{\leqslant 2}(S)$ is 
the set of interval partitions of $S$ into
\emph{at most} two  parts, and 
$\mathcal{I}_{2}(S)$ is the set of \index{partition!interval} interval partitions of $S$ into
\emph{exactly}  two  parts.\footnote{The case of interval partitions with an \emph{arbitrary} number of parts is analysed in \cite{EBMB-BaSh17}.} Clearly, 
\[
\mathcal{I}_{2}(S) = \{ \cA_k \, : \, 1 \leqslant k \leqslant n-1 \}\ts,
\]
where 
$\cA^{}_k \mathrel{\mathop:}=  \{\{1,2,\dotsc, k\},\{k+1, \dotsc, n\}\}$.
The partition $\cA^{}_k$ is the result of
a single-crossover event after site $k$. Obviously, there is a 
one-to-one correspondence between the elements of $\mathcal{I}_2(S)$
and  those of $S\setminus \{n\}$.

Likewise, there is a one-to-one correspondence between $\mathcal{I}(S)$ and the
set of subsets of $S \setminus \{n\}$. Let 
$G = \{ j^{}_1, \ldots, j^{}_{\lvert G \rvert} \} \subseteq S \setminus
\{n\}$, with $j^{}_1 < j^{}_2< \dots < j^{}_{\lvert G \rvert}$. Let then 
$\mathcal{S}(\varnothing)=\boldsymbol{1}$, and, for $G \neq \varnothing$, let
$\mathcal{S}(G)$ denote the \index{partition!interval} interval partition
\[
  \mathcal{S}(G) \,  \mathrel{\mathop:}= \, \bigl \{ \{ 1, \ldots, j^{}_1 \}, \{j^{}_1+1, \ldots, j^{}_2 \},
  \ldots, \{j^{}_{\lvert G \rvert}+1, \ldots, n \} \bigr \}\ts .
\]
In particular, ${\mathcal S} \big (S \setminus \{n\} \big )= \bigl \{\{0\},\ldots, \{n\} \bigr \}$. It is clear that 
${\mathcal S}(H) \preccurlyeq {\mathcal S}(G)$ if and only if $G \subseteq H$.
It is also obvious that $\mathcal{S}$ defines a bijection; its inverse,
\begin{equation}\label{EBMB-varphi}
\varphi = \mathcal{S}^{-1}, 
\end{equation}
associates with every \index{partition!interval} interval partition
of $S$ the corresponding subset of $S \setminus \{n\}$, so that
$\varphi(\mathcal{S}(G))=G$ for all $G \subseteq S \setminus \{n\}$.  

It is clear that $\mathcal{S}(G)$ may  alternatively be represented as
\begin{equation}\label{eq:ordered_part_repr}
  \mathcal{S}(G) \, = \, \boldsymbol{1} \wedge \cA_{j^{}_1} \wedge  \cA_{j^{}_2} \wedge \ldots \wedge \cA_{j^{}_{\lvert G \rvert}}\ts,
\end{equation}
where $\wedge$ denotes the coarsest common refinement; note that the action of $\wedge$ is commutative. In particular, 
one has $\mathcal{S} \big ( G \cup \{k\} \big ) = \mathcal{S}(G) \wedge \cA_k$. More precisely, let $\mathcal{S}(G) = 
\cB = \{B_1, \ldots, B_m \}$ and $k \in S \setminus \{n\}$. 
Then,
\[
  \cB \wedge \cA_k = \begin{cases} \cB, & k \in G\ts, \\ 
                   (\cB \setminus B_i) \cup \cA_k |^{}_{B_i}, & k \in S \setminus \big ( G \cup \{n\} \big )\ts,
       \end{cases}  
\]
where, in the latter case,
$B_i$ is the \emph{unique} block that contains $k$;  the other blocks are
not affected. 

Let us now  connect this to the \index{partitioning process} partitioning process. 
Assume that we have $\varSigma_\tau = \cB
= \{B_1, \ldots, B_m\} = \mathcal S (G)$ for some $G \subseteq S \setminus \{n\}$ and fix one index \mbox{$1 \leqslant i \leqslant m$}.
Evaluating the rates in Proposition~\ref{EBMB-detlimit} with the help of the marginal recombination rates \eqref{EBMB-marg-rates} then reveals that, in the single-crossover case, the
only (non-silent) transitions involving block $B_i$ are 
\[
\begin{split}
   & \cB  \: \mapsto  \: (\cB \setminus B_i) \cup \cA_k |^{}_{B_i}  \, = \, \cB \wedge \cA_k, \\
   & \text{at rate }   \varrho(\cA_k) \text{ for all } k \in B_i \setminus \big ( G \cup \{n\} \big )\,.
\end{split}
\]
If all blocks are taken into account, we therefore get the transitions
\begin{equation}\label{eq:scr_transitions}
\begin{split}
  & \mathcal{S}(G)  =  \cB \:  \mapsto \: \cB \wedge \cA_k = \mathcal{S} \big ( G \cup \{k\} \big ), \\
  & \text{at rate } 
  \varrho(\cA_k) \text{ for all } k \in S \setminus \big ( G \cup \{n\} \big ).
\end{split}
\end{equation}
 Since $\mathcal{S}\big ( G \cup \{k\}\big )$  is again an \index{partition!interval} interval partition, it is clear that
$\{\varSigma_\tau\}_{\tau \geqslant 0}$,  when started in
$\mathcal{I}(S)$, will never leave $\mathcal{I}(S)$. 

\begin{remark}
The property that recombination according to $\cA$ induces
the transition from $\cB$ to $\cB \wedge \cA$ is a special (and decisive)
property
of the single-crossover setting, where $\cA \in \mathcal{I}_{\leqslant 2}(S)$ and
$\cB \in \mathcal{I}(S)$, which implies that $\cA$ refines at most one block of $\cB$.
This is \emph{not} true in  the general case, where the possible refinements are considerably more complex.
\end{remark}

We are now well prepared to calculate  
$a_t$. We could work via
the matrix exponential of $Q$ and use its special structure 
resulting from the restriction to $\mathcal{I}(S)$;
however, we pursue a more elegant approach based on 
Eqs.~\eqref{eq:ordered_part_repr} and \eqref{eq:scr_transitions}.
To this end, let $\varSigma_0=\boldsymbol{1}$ and conclude 
from Eq.~\eqref{eq:scr_transitions} that
$(\varSigma_\tau)_{\tau \geqslant 0}$ is governed by the arrival
of $\cA_k$-events that happen independently of each other
at rate $\varrho(\cA_k)$. The waiting times $T_k$ until
$\cA_k$ appears are therefore independent and
exponentially distributed with parameters $\varrho(\cA_k)$. 
Let now $\cC$ be an \index{partition!interval} interval partition as in Eq.~\eqref{eq:ordered_part_repr},
that is, $\cC = \mathcal{S}(G)$ for some $G \subseteq S \setminus \{n\}$ (so $G = \varphi(\cC)$).  
Taking Eqs.~\eqref{eq:ordered_part_repr} and \eqref{eq:scr_transitions}
together, we see that $\varSigma_t = \cC$ if and only if all $\cA_k$-events with 
$k \in G$ have occurred, while all $\cA_j$-events with $j \in S \setminus \big (G \cup \{n\} \big )$
have not. We therefore get
\[
\begin{split}
   a_t^{}(\cC)  & = \, \PP(\varSigma_t = \cC \mid \varSigma_0 = \boldsymbol{1})
          \, = \, \prod_{k \in G} \PP \big  (T_k < t \big ) 
                  \prod_{\ell \in S \setminus (G \cup \{n\})} 
                \PP \big  (T_\ell \geqslant t \big ) \\
                & = \, \prod_{k \in G} \big ( 1- \ee^{-t \varrho(\cA_k)} \big )
                  \prod_{\ell \in S \setminus (G \cup \{n\})} \ee^{-t \varrho(\cA_\ell)}\ts.
\end{split}
\]
With these coefficients, Theorem~\ref{EBMB-thm-recosol} indeed turns into
an explicit and simple
solution of the recombination equation. Let us summarise our result as follows.
\begin{corollary}[Single-crossover recombination]\label{EBMB-scr}
Assume single-crossover recombination, that is, $\varrho(\cA) > 0$
implies $\cA \in \mathcal{I}_{\leqslant 2}(S)$. The probability vector \/ $a^{}_t$
from Theorem~\textnormal{\ref{EBMB-thm-recosol}}  is then given by $a^{}_t(\cC) = 0$ if $\cC \notin \mathcal{I}(S)$ and,
for $\cC \in \mathcal{I}(S)$, by
\[
  a^{}_t(\cC) \, = \, \prod_{k \in G} \big ( 1- \ee^{-t \varrho(\cA_k)} \big )
                  \prod_{\ell \in S \setminus  (G \cup \{n\}  )} \ee^{-t \varrho(\cA_\ell)}\ts,
\]
where \/ $G=\varphi(\cC)$ of Eq.~\eqref{EBMB-varphi}. 
\end{corollary}

In fact, the content of Corollary~\ref{EBMB-scr}  was originally obtained by analytic means in \cite{EBMB-BaBa03}; we have  recovered it here in genealogical terms. Note that the exponential convergence to the product measure of Corollary~\ref{EBMB-asymptotics} is obvious here from the explicit formula for the $a^{}_t(\cC)$.

 \section{Recombination  in discrete time}\label{sec:discrete}

Let us finally turn our attention to the discrete-time analogue of
Eq.~\eqref{EBMB-detreco}, namely the discrete-time \index{dynamical system} dynamical system
\begin{equation}\label{EBMB-reco-discrete}
   \omega^{}_{t+1} \, = \!  \sum_{\cA \in \cP^{}_{\! \leq 2} (S)} \!
   r (\cA) \, R^{}_{\cA} (\omega^{}_{t}) \ts,
\end{equation}
which is often considered in population genetics \cite{EBMB-Bu00, EBMB-Da02,EBMB-Ly92,EBMB-Ma17}.
Here, $t\in \NN_{0}$ now denotes
discrete time (counting generations); the  initial condition
is again \mbox{$\omega^{}_{0} \in \boldsymbol P(X)$}. The iteration describes the \emph{synchronous}
formation of a new population  from the parental one.
The parameters are now the recombination probabilities $r
(\cA)$ for  $\cA\in \cP^{}_{\! \leqslant 2} (S)$.
Obviously, $\omega^{}_{t+1}$ is a 
convex combination of $\omega^{}_t$ recombined in all possible ways,
so $\boldsymbol P (X)$ is preserved under the iteration.

In analogy with the derivation of the continuous-time recombination equation as the limit of a  finite-$N$ Moran model, the discrete-time recombination equation may be obtained as the law of large numbers of an underlying \index{Wright--Fisher model!with recombination} Wright--Fisher model with recombination; see \cite{EBMB-BaWa14} for the special case of single crossovers. Rather than working this out explicitly, we simply state the plausible fact that there is again an underlying \index{partitioning process} partitioning  process, $(
\varSigma_{\tau})_{\tau \in \NN_{0}}$. This is now a Markov chain in
\emph{discrete} time, again with values in $\cP (S)$ and starting at
$\varSigma^{}_{0} = \boldsymbol{1}$. When $\varSigma_{\tau} = \cA$, in
the time step from $\tau$ to $\tau+1$, part $A$ of $\cA$ is replaced by
$\mathfrak{a} \in \cP (A)$ with probability $r^{A}(\mathfrak{a})$,
independently for each $A \in \cA$. Note that, in contrast to the con\-tin\-u\-ous-time case, \emph{several}
parts can be refined at the same time, which makes the discrete-time case
actually more complicated. Of course, $\mathfrak{a} = \{ A\}$, which
means no action on this part, is also possible. Put together, it is
not difficult to verify that one ends up  with the Markov transition
matrix $M$ with elements
\[
    M^{}_{\cA \cB} \, = \, \begin{cases}
    \prod_{A \in \cA} r^{A} (\cB |^{}_{A}) ,
    & \text{if } \cB \preccurlyeq \cA , \\
    0 , & \text{otherwise}. \end{cases}
\]
In particular, $M = \bigl(M^{}_{\cA \cB}\bigr)^{}_{\cA,\cB\in\cP (S)}$
is a lower-triangular Markov matrix. (Let us note in passing that the triangular form, which also appears in the continuous-time case, motivated to revisit the Markov embedding problem \cite{EBMB-BaSu19}.)
The analogue of  Theorem~\ref{EBMB-thm-recosol} reads as follows.

\begin{theorem}[Solution of the discrete-time recombination equation \cite{EBMB-BaBa16}]
  \label{EBMB-thm-discreterecosol}
  The solution of the recombination equation \eqref{EBMB-reco-discrete} is given by 
\begin{equation}
     \omega^{}_{t} \, =  
     \! \sum_{\cA\in\ts\cP(S)}\! \!  a^{}_t(\cA) 
    \, R^{}_{\cA} (\omega^{}_{0}) \, = \, \EE \big ( R_{\varSigma_t}(\omega_0)   \mid \varSigma_0= \boldsymbol{1} \big ),
  \end{equation}
  where
  \[ a_t(\cA) \, = \, \PP \bigl( \varSigma_{t} = 
    \cA \mid \varSigma_0 = \boldsymbol{1} \bigr) = (M^t)^{}_{\boldsymbol{1} \cA}.  
  \] 
\end{theorem}

It is tempting to assume that, again in analogy with continuous time, the case with single crossovers might be  amenable to a simple solution. This is, however, not true. The reason is that, in continuous time, the single-crossover events appear independently by the very nature of the continuous-time process, where the probability of two events occurring simultaneously is zero.  In contrast, the single-crossover assumption in discrete time induces \emph{dependence}. Namely, a crossover between a given pair of neighbouring sites precludes a crossover between another pair of neighbouring sites in the same block. With the help of \index{M\"obius inversion} M\"obius inversion on a suitable poset of  rooted forests, a solution was obtained nevertheless, but it is of surprising complexity \cite{EBMB-BaEs18}. However, the long-term behaviour is, once more, simple:  Corollary~\ref{EBMB-asymptotics} carries over, with $\varrho$ replaced by $r$.

\subsection*{Acknowledgements}
It is our pleasure to thank Frederic Alberti for critically reading the manuscript, and two referees for helpful comments.

\end{document}